# Vehicle Routing Problem with Vector Profits (VRPVP) with Max-Min Criterion

Dongoo Lee[1] and Jaemyung Ahn[2]

*Korea Advanced Institute of Science and Technology (KAIST), 291 Daehak-Ro, Daejeon 34141, Republic of Korea*

**Abstract:** This paper introduces a new routing problem referred to as the vehicle routing problem with vector profits. Given a network composed of nodes (depot/sites) and arcs connecting the nodes, the problem determines routes that depart from the depot, visit sites to collect profits, and return to the depot. There are multiple stakeholders interested in the mission and each site is associated with a vector whose $k^{th}$ element represents the profit value for the $k^{th}$ stakeholder. The objective of the problem is to maximize the profit sum for the least satisfied stakeholder, i.e., the stakeholder with the smallest total profit value. An approach based on the linear programming relaxation and column-generation to solve this max-min type routing problem was developed. Two cases studies – the planetary surface exploration and the Rome tour cases – were presented to demonstrate the effectiveness of the proposed problem formulation and solution methodology.

**Keywords**: Vehicle Routing Problem with Profits (VRPP), Vector Profits, Multiple Stakeholders, Max-Min Criterion, Planetary Surface Exploration, Tourist Routing

## Nomenclature

| | |
|---|---|
| $n_C$ | = number of (candidate) sites for visitation |
| $n_R$ | = maximum number of routes |
| $n_S$ | = number of stakeholders |
| $i$ | = index representing a site |
| $j$ | = index representing a route (or, equivalently a set of sites) |

---

[1] PhD candidate, department of aerospace engineering

[2] Associate professor of aerospace engineering, corresponding author (Email: jaemyung.ahn@kaist.ac.kr)

| | |
|---|---|
| $k$ | = index representing a stakeholder |
| **C** | = set of candidate customers |
| **J** / **J**$_f$ | = index set of all routes / feasible routes |
| **R**$_j$ | = set of exploration sites belonging to route $j$ |
| TSP$_j$ | = length/time to complete route $j$, found by solving the traveling salesman problem |
| $p_i^k$ | = profit obtainable at site $i$ for stakeholder $k$ |
| $t_i$ | = time required to obtain the profit at site $i$ |
| **r**$_i$ ($x_i$, $y_i$) | = locations of base ($i = 0$) and customers ($i = 1, \ldots, n_C$) |
| **b**$_r$ | = resource consumption budget (limit) *for a single-route* |
| **b**$_m$ | = resource consumption budget (limit) *for the whole mission* |
| **c**$_d$ / **c**$_r$ | = on-arc / on-site resource consumption coefficient vectors |
| **d**$_d$ / **d**$_r$ | = on-arc/on-site resource consumption coefficient vectors *for the whole mission* |
| **h**$_j$ | = resource consumption associated with route $j$ |
| $x_j$ | = binary decision variable; equals 1 if the route $j$ is included in the solution and 0 otherwise |
| **x** | = vector of $x_j$ values |
| $g_o$ | = optimality gap of the near-optimal solution |
| $J$ | = objective function for the optimization (= minimum sum of profits for each profit type) |
| * | = superscript representing the optimal value |
| $f$ | = subscript representing a feasible route |

# I.     Introduction

The vehicle routing problem (VRP) has been actively studied and used to address operational challenges, primarily in the field of supply chain management (SCM), since it was originally introduced by Dantzig and Ramser (1959) as an extension of a classical traveling salesman problem (TSP). Given a network composed of nodes (a depot and sites) and arcs between them, the original VRP formulation determines multiple routes of a vehicle (a delivery truck) that 1) depart from and return to the depot and 2) collectively visit all sites. The objective of the problem is to minimize the total travel distance under constraints such as the delivery capacity of the vehicle and the number of routes.

The vehicle routing problem with profits (VRPP) is an important variant of VRP. In the VRPP, 1) a numerical "profit value" representing the amount of benefit obtained by the visit is assigned to each site, and 2) leaving certain (less attractive) sites unvisited is allowed. Therefore, instead of minimizing the total travel distance over the routes that collectively visit all sites, the VRPP maximizes the sum of collected profit values by determining the sites to visit (a subset of the whole sites) and their visiting routes while satisfying the constraints such as the resource consumption on each route or the overall mission.

In the original VRPP formulation, the profit obtainable by visiting a site was defined as a *scalar* value. This definition implies that there exists a *single* stakeholder interested in the routing mission and the stakeholder can determine the profit function relating the site visit to a numerical value. In reality, however, some routing problems involve multiple stakeholders whose profit functions are quite different. One example is the determination of touring routes for a group composed of individuals with different interests (e.g. history, architecture, and food). In this case, the profit assigned to a site should be a *vector* whose elements represent the benefits to individual members of the group.

This paper introduces a new routing problem referred to as the *vehicle routing problem with vector profits* (VRPVP) that can reflect the perspectives of multiple stakeholders by introducing the concept of vector profits, which is its key contribution. The objective of the problem is to maximize the profit sum of the least satisfied stakeholder. A max-min binary programming procedure composed of the linear programming (LP) relaxation and the column generation technique is proposed to find a near-optimal solution of the problem and its optimality gap. Numerical experiments and two cases studies are conducted to demonstrate the effectiveness of the proposed problem formulation and solution methodology.

The rest of this paper is organized as follows. Section II provides the review of the past studies related to the subject of this paper. Section III presents the mathematical formulation of the VRPVP by defining its objective function and constraints. Section IV explains the steps of a procedure to obtain the near-optimal solution of the VRPVP based on a column-generation technique. The effectiveness of the formulation and the solution procedure are validated through numerical experiments (Section V) and two case studies (Section VI) – the planetary surface exploration and the Rome tour cases. Finally, Section VII presents the conclusions of the study and discusses potential future research.

# II. Literature Review

Many variations and extensions of the original VRP have been developed to address real-world applications effectively. Consideration of multiple depots (multi-depot VRP, MDVRP), integrated decisions on routing and depot selection (location routing problem, LRP), and introduction of new constraints on site visit time windows (VRP with time windows, VRPTW) are examples of these variations (Laporte et al., 1988; Laporte et al., 1989; Nagy and Salhi, 2007; Bräysy and Gendreau, 2005; Berger et al. 2007). Toth and Vigo (2014) provided a very comprehensive survey on these variations in their book.

A number of studies on applications and solution methodologies for *the routing problem with profits*, which is one of these variants, can be found in literature. For example, Balas (1989) introduced the prize-collecting traveling salesman problem (PCTSP) and discussed the methods to obtain its exact solution. Chao et al. (1996) formulated the routing of multiple agents/members to maximize the sum of scores obtainable by visiting sites (the team orienteering problem, TOP) and solved the problem using a heuristic algorithm. Butt and Cavalier (1994) and Butt and Ryan (1999) provided procedures to solve the VRPP using a heuristic algorithm and a column generation based optimization technique, respectively.[3] Chu (2005) introduced a routing problem for two different carrier types (truckload and less-than-truckload) reflecting the viewpoint of a logistics manager responsible for the decisions on the carrier type and routing to minimize the total cost, and proposed a heuristic algorithm to solve the problem. Another interesting example is the tourist trip design problem (TTDP) introduced by Vansteenwegen and Van-Oudheusden (2007), which can be used as a real-time algorithm for mobile applications. The TTDP is an extension of the TOP that can consider the factors that are important for traveler such as time windows, budget limitations, attraction values, and scenic routes. A survey on various approaches to find the solution of the VRPP such as exact methods, classical heuristic procedures, and metaheuristics was provided by Feillet et al. (2005).

Some relatively recent studies on vehicle routing address the cases with multiple objectives – instead of a single objective adopted in traditional problems (e.g. maximizing total travel distance and minimizing total cost). Jozefowiez et al. (2008) proposed a meta-heuristic algorithm to obtain the Pareto optimal solutions of a bi-objective traveling salesman problem with profits (TSPP), whose two objectives were minimizing the total tour length and maximizing the total profits. Schilde et al. (2009) introduced a bi-objective orienteering problem considering two different benefit categories obtainable by visiting a site. Two metaheuristics-based algorithms – the Pareto ant colony optimization (P-ACO) and the Pareto variable neighborhood search (P-VNS) – to generate the non-dominated front of the problem were proposed. As the most recent study, Matl et al. (2017) formulated the personal scheduling for a mobile freelancer as the bi-objective orienteering problem. Two objectives considered in their study were to maximize the task planning and to maximize the enjoyable free time. Non-dominated solutions of the problem were obtained by a metaheuristic algorithm based on large neighborhood search (LNS). Note that all the aforementioned multi-objective studies addressed two objectives – the authors could not find any routing problems that deal with three or more objectives.

---

[3] The problem was referred to as the multiple tour maximum collection problem (MTMCP) in their original paper.

The problem proposed in this paper is similar to the robust VRP in that the concept of profit vector relates multiple values to a site. The robust VRP can address the uncertainty in the parameters and data of routing problems (e.g. demand, travel cost, and service time), which are treated as deterministic in traditional VRPs. Sungur et al. (2008) presented a robust vehicle routing problem to minimize the cost while satisfying the customers' demand under uncertainty. Ordonez (2010) and Solano-Charris (2015) provided comprehensive surveys on various uncertainty models (such as costs, demand, time, and customers) and solution methodologies for the robust VRP in their paper, respectively.

The max-min criterion has been used in some studies on robust VRP – conceptually or explicitly, while the authors could not find any VRPP study that adopted this criterion. Han et al. (2013) considered multiple scenarios associated with various forecasts (on uncertain travel time) for route selection to improve the performance of the routing mission under the worst-case scenario. Ogryczak (1997) applied the lexicographic minimax approach to a location selection problem. He pointed out that the approach can overcome the criticism on the traditional minimax problem that the efficiency (Pareto optimality) of the solution is not guaranteed.

# III. Problem Description

## A. Vehicle Routing Problem with Vector Profits (VRPVP)

As was mentioned in Section I of this paper, the VRPVP is a variant of the vehicle routing problem with profits (VRPP). Given the depot/sites locations and profit values associated with the sites, the VRPP determines the set of routes starting and terminating at the depot that maximize the sum of profits obtained from the sites visited by the routes under certain constraints. Each site can be visited at most once and it is not necessary to visit all sites, which differentiates the VRPP from classical VRPs, which require visits to all sites. Figure 1-(a) illustrates a VRPP instance and one example of a feasible (but not necessarily optimal) solution. For each site (circle), a profit value (number printed near the circle) is assigned, and the sum of the profit values obtained from the visited sites is calculated to obtain the objective function (total profit sum) of the VRPP.

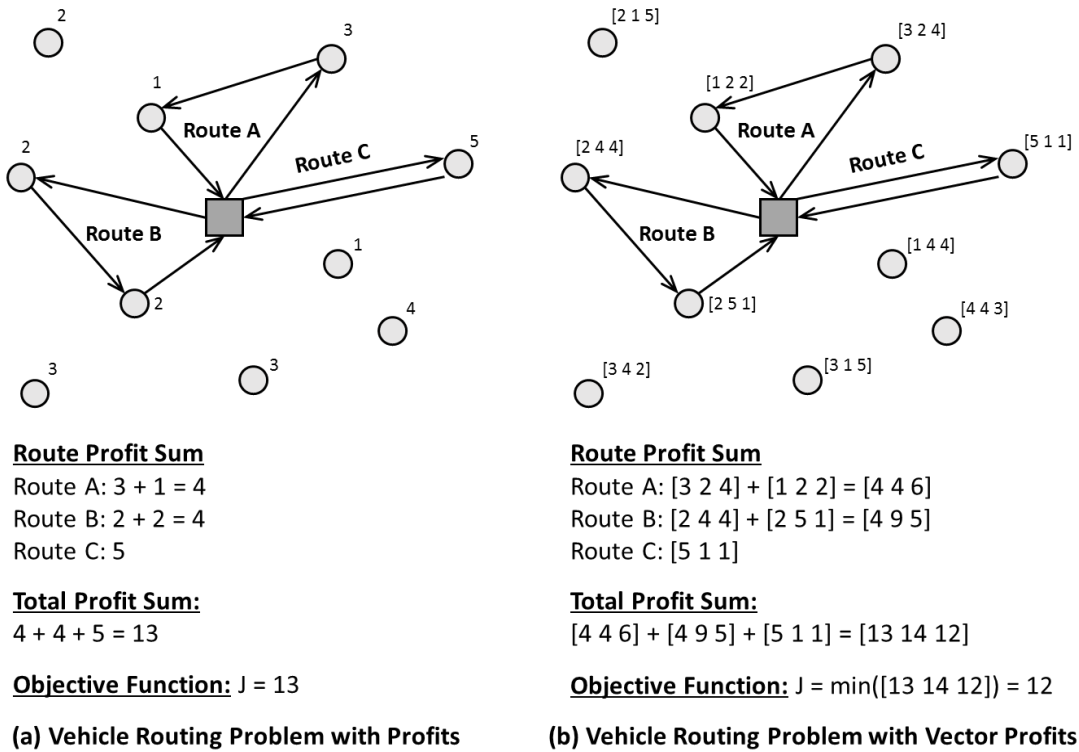

Figure 1: Sample problem instances – (a) VRPP, (b) VRPVP

In certain applications involving multiple stakeholders, however, representing the profit of visiting a site as a scalar value is not sufficient to represent the perspectives of different stakeholders on the site. To address this challenge, we propose a new problem referred to as the vehicle routing problem with vector profits (VRPVP), which uses a *vector* whose dimension is equal to the number of

stakeholders – to express the profit obtainable by visiting a site. In this expression, component $k$ of the vector represents the profit value obtained by stakeholder $k$. Fig. 1-(b) presents a sample instance of the VRPVP. A three-dimensional vector is assigned to each site, which indicates that there are three different stakeholders involved in the routing problem. For example, the vector [2 1 5] is assigned to a site located in upper left corner, which means that the first stakeholder obtains the profit value of 2, the second obtains 1, and the third obtains 5, by visiting the site.

Since the profit assigned to a site becomes a vector, the sum of the profits for a route or for the overall mission is a vector whose $k^{th}$ component is the sum of the profits for the $k^{th}$ stakeholder. In Fig. 1-(b), for example, the total profit sum obtained from the mission (composed of routes A, B, and C) is [13 14 12]. This finding means the three stakeholders obtained the total profit values of 13, 14, and 12, respectively, from the mission.

With the total profit sum in vector form, how to formulate the VRPVP mathematically as an optimization problem – definition of its objective function, in particular – remains an issue. Several different approaches for this issue can be considered. One approach is to calculate the (weighted) sum of the elements of the total profit sum vector and define it as an objective of the routing problem. The weighted sum approach can transform the VRPVP into traditional VRPP and solve by existing solution methodologies, which is its key advantage. However, the approach has a critical limitation that low weighted elements are likely to be ignored during the decision making procedure – like the drawback of the majority rule. To define a multi-objective optimization (MOO) problem whose objectives are the elements of the vector – total profit sums for different stakeholders – is another approach. This approach generates the family of efficient solution referred to as the Pareto front, which can reflect the perspectives of all stakeholders. However, generating the Pareto front for multi-objective routing problem is very difficult and may be not tractable sometimes.

The objective function used in the VRPVP formulation presented in this paper is the maximization of the minimum element of the total profit sum vector (the Max-Min criterion). The procedure to calculate the objective function is described in Fig. 1-(b). Out of the elements of the vector [13 14 12], the total profit sum for stakeholder 3 (= 12) is the minimum, and is determined to be the objective function (to maximize) of the problem.

## B. Mathematical Formulation of the Problem

This subsection provides the mathematical formulation of the VRPVP. To define a routing problem, we assume that a depot and multiple sites ($C = \{1, \ldots, n_c\}$) are given. For each site $i$, a profit vector $\mathbf{p}_i$ is defined as follows.

$$\mathbf{p}_i = [p_i^1, \cdots, p_i^{n_s}], \qquad (1)$$

where $p_i^k$ denotes the profit value obtainable in site $i$ for stakeholder $k$ and $n_s$ is the number of stakeholders. Note that we are interested in the cases with multiple stakeholders, and $n_s$ is an integer no smaller than 2.

A vehicle visits the subset of $C$ under resource constraints (e.g. amount of fuel consumed, time spent, and number of routes) imposed on both individual routes and the whole mission. Note that the resource consumption can be classified into two types: *on-arc* and *on-site* consumption types. The on-arc consumption arises while the vehicle is moving from one site to another, and is directly proportional to the *length of the route* or the *time spent on the route*. The on-site consumption – proportional to the *time spent on site* – takes place when the profit from the site is acquired.

We additionally assume that only the shortest closed path that visits the given subset of sites and depot locations, which is the solution of the traveling salesman problem (TSP), is selected. Once the sites to be visited are specified, a path that starts at the depot, visits all the specified sites, and returns to the depot is determined as an associated route. In this regard, we can define the index set of possible routes ($\mathbf{J}$) as follows.

$$\mathbf{J} = \{0, 1, \cdots, 2^{n_C} - 1\} \qquad (2)$$

An integer $j \in \mathbf{J}$ is associated with a route such that site $i$ belongs to the route if the $i^{th}$ digit ($d_i^j \equiv \mathrm{mod}([j/2^{i-1}], 2)$) is equal to 1, and does not otherwise.[4] With this information, we can define $\mathbf{R}_j$ as the set of customers that belongs to route $j$.

In addition, we define the index set of feasible routes $\mathbf{J}_f$ with in-route consumption coefficient vectors $\mathbf{c}_d$ and $\mathbf{c}_r$ and the budget vector $\mathbf{b}_r$ as follows:

---

[4] $[x]$ is the largest integer that does not exceed $x$ and $\mathrm{mod}(a, m)$ denotes the remainder after division of $a$ by $m$.

$$\mathbf{J}_f = \{ j \in \mathbf{J} \mid \text{TSP}_j \cdot \mathbf{c}_d + (\sum_{i \in \mathbf{R}_j} t_i) \cdot \mathbf{c}_r \leq \mathbf{b}_r \} \tag{3}$$

where TSP$_j$ is the solution of the traveling salesman problem with sites in $\mathbf{R}_j$ and the depot, and $t_i$ is the time required to stay at site $i$ to obtain the profit. Any types of resources required to operate the vehicle can be considered for constraints – these are the elements of $\mathbf{b}_r$.

Note that the TSP$_j$ denotes the *shortest travel length* or the *minimum travel time* to complete the routes depending on the type of cost (distance or time) associated with the arc representing the movement between two sites. In addition, $\mathbf{c}_d$ and $\mathbf{c}_r$ respectively represent *route-constraining* resource consumption per unit length (on-arc) and per unit time (on-site). For example, if the total time to complete the given mission (in hours) is considered as the constraining resource, and TSP$_j$ is the shortest travel length (in km), then $\mathbf{c}_d$ is the inverse of the vehicle's velocity (in hour/km) and $\mathbf{c}_r$ is unity (in hour/hour).

Mathematical formulation of the VRPVP as a max-min type integer-programming problem is presented as follows.

$$[\mathcal{P}_0: \text{VRPVP}] \quad J_O = \max_{x_j} \min_k \sum_{j \in \mathbf{J}_f} \left( r_j^k \cdot x_j \right) = \max_{x_j} \min_k \sum_{j \in \mathbf{J}_f} \left( \left( \sum_{i \in \mathbf{R}_j} p_i^k \right) \cdot x_j \right), \tag{4}$$

subject to

$$\sum_{j \in \mathbf{J}_f} (\mathbf{a}_j \cdot x_j) \leq \mathbf{1}_{n_E}, \tag{5}$$

$$\sum_{j \in \mathbf{J}_f} (\mathbf{h}_j \cdot x_j) \leq \mathbf{b}_m, \tag{6}$$

$$\sum_{j \in \mathbf{J}_f} x_j \leq n_R, \tag{7}$$

$$x_j \in \{0, 1\}. \tag{8}$$

The objective function of the VRPVP, which is maximization of the minimum total profit sum obtained from the mission for all stakeholders, is expressed in Eq. (4). Decision variables for this problem are $x_j$ and $k$. The binary variable $x_j$ is equal to 1 if route $j$ – the TSP solution for the depot and sites belonging to $\mathbf{R}_j$ – is included in the solution, and 0 otherwise. The variable $k$ is the index of the stakeholder that has the minimum total profit sum. Also, $r_j^k$ is the sum of profits for stakeholder $k$ obtained from route $j$ (i.e., obtained from sites belonging to $\mathbf{R}_j$). Eq. (5) imposes the constraint that a site cannot be visited more than once. In this equation, $\mathbf{a}_j$ is the $n_C$-dimensional column vector whose $i^{th}$ element is 1 if $i \in \mathbf{R}_j$ and 0 otherwise, and $\mathbf{1}_{n_C}$ is the $n_C$-dimensional column vector whose elements are all 1's. Eq. (6) expresses resource constraints over the whole mission. The vector $\mathbf{h}_j$ represents the amount of resource consumed while the vehicle is traveling on route $j$, and is defined as follows.

$$\mathbf{h}_j = \text{TSP}_j \cdot \mathbf{d}_d + (\sum_{i \in \mathbf{R}_j} t_i) \cdot \mathbf{d}_r, \qquad (9)$$

where $\mathbf{d}_d$ and $\mathbf{d}_r$ respectively represent *mission-constraining* resource consumption per unit length (on-arc) and per unit time (on-site). We can obtain the total consumption of the mission-constraining resource by summing $\mathbf{h}_j$ over all $j$ included in the solution, which should be no greater than the budget for the mission-constraining resource ($\mathbf{b}_m$). Note that we can set up the mission-constraining resources (elements of $\mathbf{b}_m$) differently from the route-constraining resources (elements of $\mathbf{b}_r$), which makes $\mathbf{c}_d/\mathbf{c}_r$ different from $\mathbf{d}_d/\mathbf{d}_r$. For example, it is possible that a route constrains the amount of fuel used to complete the route but the mission constrains the total time spent on the whole mission. Eq. (7) expresses the constraint that the number of routes for the mission should be no greater than the maximum value ($n_R$), and Eq. (8) represents the constraint that $x_j$ is binary.

The max-min type formulation for the VRPVP described in Eqs. (4)-(8) can be rewritten in a simple minimization mixed-integer linear programming (MILP) formulation, which will be used for its solution methodology. The MILP formulation of the VRPVP ($\mathcal{P}$) is as follows.

[$\mathcal{P}$: VRPVP – MILP Formulation] $\qquad J_O = \min_{x_j, z} (-z), \qquad (10)$

subject to Eqs. (5)-(8) and an additional constraint,

$$z \leq \sum_{j \in \mathbf{J}_f} \left( r_j^k \cdot x_j \right) \quad (1 \leq k \leq n_S). \qquad (11)$$

In this formulation, we introduce a new real-value decision variable $z$ to convert the max-min formulation to a minimization problem (Bertimas and Tsitsiklis 1997, pp. 16-17). The objective of the new problem is to minimize -*z* (or to maximize *z*) where *z* is no greater than the total profit sum of any stakeholder (Eq. (11)). The solution methodology presented in the next subsection is developed based on this MILP formulation.

## IV. Solution Methodology for the VRPVP

We developed the procedure to solve the problem by modifying the methodology used to solve the VRPP by Ahn (2008). The procedure is designed based on the column-generation technique (Bertsimas and Tsitsiklis 1997; Simchi-Levi *et al*. 2005; Chabrier 2006). The technique is widely used to solve linear programming (LP) problems for which the system matrix has too many columns, making the full enumeration prohibitively difficult. The property that only very few columns among the large number of possible columns are generally included in the final solution is utilized in this technique. Feillet et al. (2005) provided a good summary on how to apply the column generation technique to vehicle routing problems. The number of columns of the VRPVP (problem $\mathcal{P}$) is the same as the cardinality of $\mathbf{J}_f$, the maximum value of which is $2^{n_c}$, and the column-generation method can be effectively utilized.

The outline of the procedure is explained as follows. First, the LP relaxation of the original MILP, $\mathcal{P}_L$, is obtained and solved to optimality by using the column-generation technique. The columns generated to find the optimal solution of $\mathcal{P}_L$ are stored and used to construct the fractional MILP of the original problem, $\mathcal{P}_A$. The fractional MILP was solved to optimality by using the branch and bound technique, and this solution is the near-optimal solution of the original VRPVP. The worst-case optimality gap between this near-optimal solution and the true optimum is computed by comparing the optimal objective functions of $\mathcal{P}_L$ and $\mathcal{P}_A$. A block diagram describing the proposed solution procedure is presented in Fig. 2. The rest of this subsection discusses the stepwise details of the solution methodology from Step 2[5].

---

[5] Step 1 of the procedure has been explained in the previous subsection (III-B).

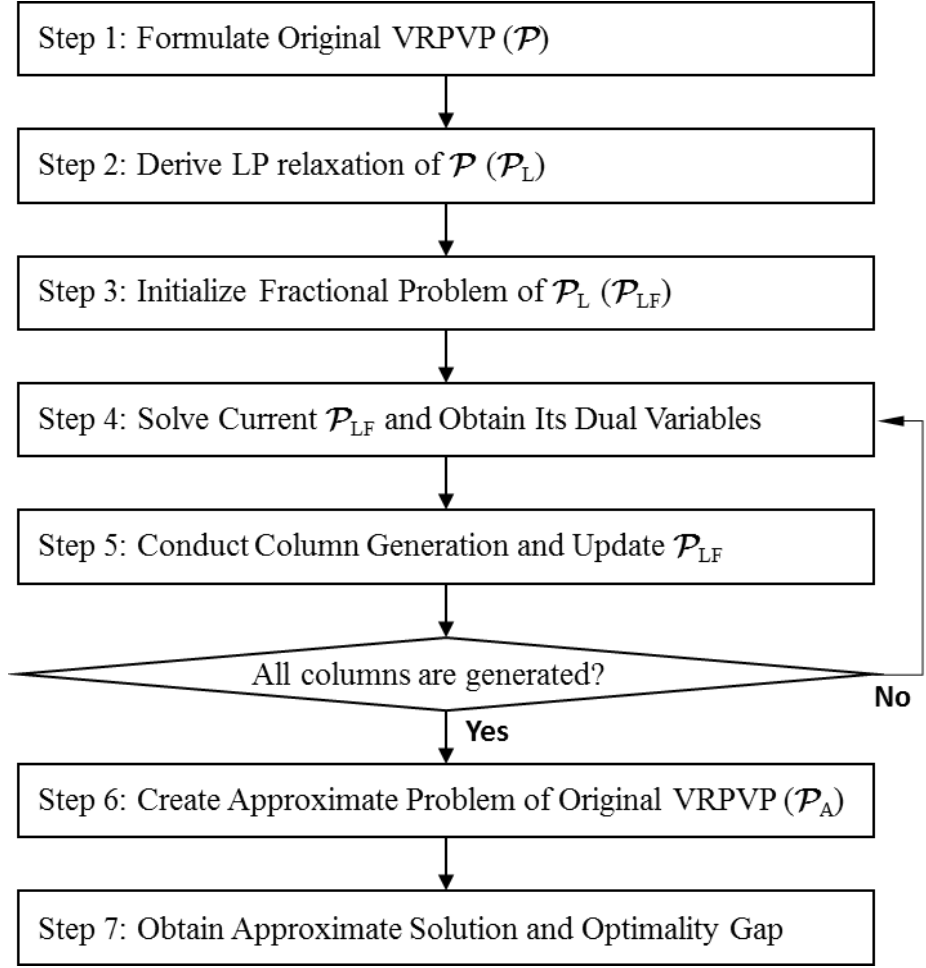

**Figure 2: Block diagram of the solution procedure for VRPVP**

**(Step 2) Derive the LP relaxation of $\mathcal{P}$ ($\mathcal{P}_L$)**

The step starts by introducing the mathematical formulation of the LP relaxation of the VRPVP ($\mathcal{P}_L$).

[$\mathcal{P}_L$: LP Relaxation of $\mathcal{P}$] $\qquad J_L = \min_{\mathbf{x},z} (-z)$, $\qquad\qquad(12)$

subject to

$$\mathbf{A}\mathbf{x} \leq \mathbf{1}_{n_C}, \qquad(13)$$

$$\mathbf{H}\mathbf{x} \leq \mathbf{b}_m, \qquad(14)$$

$$\mathbf{1}_{\mathbf{x}}^T \mathbf{x} \leq n_R, \qquad(15)$$

$$-\mathbf{R}\mathbf{x}+\mathbf{1}_{n_s} z \leq \mathbf{0}_{n_s}, \tag{16}$$

$$\mathbf{x} \geq 0. \tag{17}$$

In this formulation, the vector $\mathbf{x}$ and matrices $\mathbf{A}$, $\mathbf{H}$, $\mathbf{R}$ are defined so that Eqs. (12)-(16) are equivalent to the objective function and constraints $\mathcal{P}$ as follows.

$$\mathbf{x} \equiv [\cdots x_j \cdots]^T \quad (j \in \mathbf{J}_f), \tag{18}$$

$$\mathbf{A} \equiv [\cdots \mathbf{a}_j \cdots] \quad (j \in \mathbf{J}_f), \tag{19}$$

$$\mathbf{H} \equiv [\cdots \mathbf{h}_j \cdots] \quad (j \in \mathbf{J}_f), \tag{20}$$

$$\mathbf{R} \equiv \begin{bmatrix} \cdots & r_j^1 & \cdots \\ \cdots & \vdots & \cdots \\ \cdots & r_j^k & \cdots \end{bmatrix} \quad (j \in \mathbf{J}_f). \tag{21}$$

In addition, $\mathbf{1}_\mathbf{x}^T$ is a row vector of ones that has the same length as $\mathbf{x}$. The binary constraint of $\mathcal{P}$ (Eq. (8)) is relaxed to a nonnegative constraint expressed as Eq. (17).[6]

## (Step 3) Initialize a Fractional Problem of $\mathcal{P}_L$ ($\mathcal{P}_{LF}$)

Note that each column of $\mathcal{P}$ or $\mathcal{P}_L$ is associated with a route – indexed by $j$. Although the total number of columns of the problem is very large (= $\|\mathbf{J}_f\|$), only a very small fraction of them are relevant to the optimal solution. Using this property, we construct a fractional problem of $\mathcal{P}_L$ with initial columns that can be included in the feasible solution (e.g. simple round trips from the depot to each site), and systematically identify the candidate columns that can be included in the optimal solution of the problem; this procedure is a column-generation method. Problem $\mathcal{P}_{LF}$, the fractional problem of $\mathcal{P}_L$, is defined as follows.

---

[6] A direct relaxation of the binary constraint is $\mathbf{0} \leq \mathbf{x} \leq \mathbf{1}$. In this case, however, Eq. (13) automatically guarantees that $\mathbf{x}$ is no greater than $\mathbf{1}$, and the constraint $\mathbf{x} \geq \mathbf{1}$ can be omitted.

[$\mathcal{P}_{\text{LF}}$: Fractional Problem of $\mathcal{P}_\text{L}$]  $\quad J_{LF} = \min_{\mathbf{x},z} (-z)$,  (22)

subject to

$$\mathbf{A}_c \mathbf{x} \leq \mathbf{1}_{n_C},$$  (23)

$$\mathbf{H}_c \mathbf{x} \leq \mathbf{b}_m,$$  (24)

$$\mathbf{1}_\mathbf{x}^T \mathbf{x} \leq n_R,$$  (25)

$$-\mathbf{R}_c \mathbf{x} + \mathbf{1}_{n_s} z \leq \mathbf{0}_{n_s},$$  (26)

$$\mathbf{x} \geq \mathbf{0},$$  (27)

where $\mathbf{A}_c$, $\mathbf{H}_c$, and $\mathbf{R}_c$ are fractional matrices and $\mathbf{1}_c^T$ is a fractional row vector associated with the generated columns.

## (Step 4) Solve the Current $\mathcal{P}_{\text{LF}}$ and Obtain Its Dual Variables

The column-generation procedure identifies the columns to update $\mathcal{P}_{\text{LF}}$ so that the true optimum of problem $\mathcal{P}_\text{L}$ can be obtained by solving problem $\mathcal{P}_{\text{LF}}$. The procedure utilizes the property that the optimality of the primal problem ($\mathcal{P}_\text{L}$) indicates the feasibility of its dual ($\mathcal{P}_\text{D}$) (Bertimas and Tsitsiklis 1997), which is defined as follows.

[$\mathcal{P}_\text{D}$: Dual Problem of $\mathcal{P}_\text{L}$]  $\quad J_D = \max_{\mathbf{q}_1,\mathbf{q}_2,q_3,\mathbf{w}} (\mathbf{q}_1^T \mathbf{1}_{n_C} + \mathbf{q}_2^T \mathbf{b}_c + q_3 n_R + \mathbf{w}^T \mathbf{0}_{n_S})$,  (28)

subject to,

$$\mathbf{q}_1^T \mathbf{A} + \mathbf{q}_2^T \mathbf{H} + q_3^T \mathbf{1}_\mathbf{x} - \mathbf{w}^T \mathbf{R} \leq \mathbf{0}_\mathbf{x}^T,$$  (29)

$$\mathbf{w}^T \mathbf{1}_{n_S} = -1,$$  (30)

$$\mathbf{q}_1, \mathbf{q}_2, q_3, \mathbf{w} \leq 0.$$  (31)

In this formulation, $\mathbf{q}_1$, $\mathbf{q}_2$, $q_3$, and $\mathbf{w}$ are dual variables associated with constraints (13), (14), (15), and (16), respectively.

**(Step 5) Conduct Column Generation and Update $\mathcal{P}_{LF}$**

If the optimal solution of $\mathcal{P}_{LF}$ composed of columns that have been generated to date is the true optimum of $\mathcal{P}_L$, its dual variables ($\mathbf{q}_{1c}^*$, $\mathbf{q}_{2c}^*$, $q_{3c}^*$, and $\mathbf{w}_c^*$) should satisfy the constraint presented as Eq. (29). The existence of any column (or, equivalently, routes $j$) violating Eq. (29) for the dual variables indicates that it should be included in $\mathcal{P}_{LF}$ for further iteration. Hence, we can identify such columns (or indices $j$) using the following inequality in the column-generation procedure.

$$(\mathbf{q}_{1c}^*)^T \mathbf{A}_j + (\mathbf{q}_{2c}^*)^T \mathbf{H}_j + q_{3c}^* - (\mathbf{w}_c^*)^T \mathbf{R}_j > 0, \tag{32}$$

where $\mathbf{A}_j$, $\mathbf{H}_j$, and $\mathbf{R}_j$ are the columns of the matrices $\mathbf{A}$, $\mathbf{H}$, and $\mathbf{R}$ associated with route $j$, respectively. Eq. (32) can be reorganized using the following relationships.

$$(\mathbf{q}_{1c}^*)^T \mathbf{A}_j = (\mathbf{q}_{1c}^*)^T \mathbf{a}_j = \sum_{i \in \mathbf{R}_j} q_{1ci}, \tag{33}$$

$$(\mathbf{q}_{2c}^*)^T \mathbf{H}_j = (\mathbf{q}_{2c}^*)^T \mathbf{h}_j, \tag{34}$$

$$-(\mathbf{w}_c^*)^T \mathbf{R}_j = -\sum_{k=1}^{n_S} w_{c,k}^* \cdot r_j^k = -\sum_{i \in \mathbf{R}_j} \sum_{k=1}^{n_S} (w_{c,k}^* \cdot p_i^k), \tag{35}$$

where $q_{1ci}$ denotes the $i^{th}$ element of $\mathbf{q}_{1c}$ and $w_{c,k}^*$ is the $k^{th}$ element of $\mathbf{w}_c^*$. Then the column-generation condition becomes

$$\sum_{i \in \mathbf{R}_j} (q_{1ci}^* - \sum_{k=1}^{n_S} w_{c,k}^* p_i^k) + (\mathbf{q}_{2c}^*)^T \mathbf{h}_j + q_{3c}^* > 0. \tag{36}$$

Since $\mathbf{q}_{2c}^*$ and $q_{3c}^*$ are non-positive from Eq. (31) and the resource consumption vector $\mathbf{h}_j$ is positive, Eq. (36) can be true only if $\sum_{i \in \mathbf{R}_j} (q_{1ci}^* - \sum_{k=1}^{n_S} w_{c,k}^* p_i^k)$ is positive; this property is used for generating new columns. First, with the current optimal dual value of $\mathcal{P}_{LF}$, the values of $(q_{1ci}^* - \sum_{k=1}^{n_S} w_{c,k}^* p_i^k) \equiv u_i$ for all sites (equivalently, for all indices $i \leq n_S$) are calculated and sorted in a descending order. The sites with positive $u_i$ values are combined to identify candidate routes. To identify all feasible routes systematically, the sites are combined in a lexicographical way to construct routes (Ahn 2008, Ahn et al. 2012), and the constructed routes are tested whether they belong to $\mathbf{J}_f$ using the single-route resource constraints expressed in Eq. (3). Columns associated with the feasible routes (or *j*) are generated to update the fractional LP problem ($\mathcal{P}_{LF}$). This column-generation procedure continues until there is no additional feasible route – with respect to the single-route resource constraint – violating the dual feasibility condition described by Eq. (36).

**(Step 6) Create an Approximate Problem of the Original VRPVP ($\mathcal{P}_A$)**

The columns for the final fractional problem and associated matrices ($\mathbf{A}_f$, $\mathbf{H}_f$, and $\mathbf{R}_f$) are used to construct a MILP that can provide a near-optimal solution of the original VRPVP ($\mathcal{P}$) as follows.

[$\mathcal{P}_A$: Approximate formulation of $\mathcal{P}$] $\qquad J_A = \min_{\mathbf{x},z}(-z),$ (37)

subject to

$$\mathbf{A}_f \mathbf{x} \leq \mathbf{1}_{n_S}, \qquad (38)$$

$$\mathbf{H}_f \mathbf{x} \leq \mathbf{b}_c, \qquad (39)$$

$$\mathbf{1}_\mathbf{x}^T \mathbf{x} \leq n_R, \qquad (40)$$

$$-\mathbf{R}_f \mathbf{x} + \mathbf{1}_{n_s} z \leq \mathbf{0}_{n_s}, \qquad (41)$$

$$\mathbf{x}: \text{binary}. \qquad (42)$$

The solution of $\mathcal{P}_A$ can be obtained using one of the standard approaches for LP problems involving integer variables (e.g. the branch and bound method).

**(Step 7) Obtain the Approximate Solution and Optimality Gap**

While the optimal solution of $\mathcal{P}_{LF}$ with the final columns is same as that of $\mathcal{P}_L$, it is not guaranteed that the optimal solution of $\mathcal{P}_A$ is optimal for the original VRPVP formulation ($\mathcal{P}$) as well. Using the solution of $\mathcal{P}_{LF}$, however, we can compute a metric related to the optimality of the approximate solution obtained by the proposed procedure.

Let $J^*$ be the true optimal objective function of the original VRPVP ($\mathcal{P}$). Also, let $J_A^*$ and $J_L^*$ be the optimal objective function values obtained by solving problems $\mathcal{P}_A$ and $\mathcal{P}_L$ (or, equivalently, the final $\mathcal{P}_{LF}$ obtained after completing the column-generation procedure). The relationship between these three values can be expressed as the following inequalities.

$$J_A^* \leq J^* \leq J_L^* \tag{43}$$

The first inequality holds because $\mathcal{P}_A$ is constructed using a fraction of the columns of $\mathcal{P}$, and the second inequality holds because $\mathcal{P}_L$ is a relaxation of $\mathcal{P}$ and hence has the better objective function. These inequalities are used to compute the optimality gap between the approximate solution and the true optimal solution as follows.

$$G_O = \left(\frac{J_L^* - J_A^*}{J_L^*}\right) \times 100 \ (\%). \tag{44}$$

Note that the optimality gap expressed as Eq. (44) provides the *worst-case* gap between the approximation and the true optimum.

The next two sections present the results of the numerical experiments and two case studies demonstrating the effectiveness of the problem formulation and solution procedure introduced in this subsection.

# V.  Numerical Experiments

This section presents the results of numerical experiments for demonstrating the validity of the mathematical formulation and solution procedure for the VRPVP proposed in this paper. The benchmarking Team Orienteering Problem (TOP) instances presented in Chao et al. (1996) were used for the experiments. Note that the original instances do not involve vector profits – each site is assigned a scalar profit value. The instances were converted to VRPVP instances by 1) introducing three more stakeholders (total four stakeholders), and 2) assigning profits to sites (associated with a specific stakeholder) by shuffling the original profit values. It is guaranteed that the sums of total profits associated with different stakeholders are all the same. Instances with different positions of departure/arrival nodes were modified by co-locating the two nodes. Other parameters such as number of sites, number of routes, and resource budget were kept the same as original.

The core of the VRPVP formulation and the column-generation procedure were implemented in C. The optimal solutions of the linear fraction problems ($\mathcal{P}_{LF}$) and the approximate MILP ($\mathcal{P}_A$) were obtained using IBM ILOG CPLEX Optimization Studio. The computations were conducted on an Intel i5 quad-core CPU (2.70 GHz) with 16 GB RAM under the Windows 7 operating system. A subroutine to obtain the exact solution of traveling salesman problems developed by Carpento et al. (1995) – based on *assignment problem relaxation* and *subtour elimination branching scheme* – were used to calculate the TSP solution for visiting the site set $j$ (TSP$_j$).

Table 1 summarizes the results of the numerical experiments. Total 45 instances created with 7 different networks and various values of resource consumption budget ($b_r$) were used for the experiments. The numbers of sites of the instances range between 19 and 100. Note that only one resource type (travel distance) is considered in the test cases and $b_r$ is set as a scalar.

It can be observed that the proposed algorithm can obtain solutions for all the instances – including the one with relatively large number of sites (up to 100). The magnitudes of optimality gaps were very small; we had 0 % gap for 33 out of 45 instances and its maximum value was 3.05 %. The computing times become larger as the value of resource consumption budget ($b_r$) increases, which entails the expansion of the solution space and requires more column generations. In-depth analysis on the relationships between the characteristics of a problem instance and the optimality gap / computing time can be a potential subject for future study.

**Table 1: Results of Numerical Experiments for Benchmarking Problems**

| Instance (Chao et al., 1996) | Number of Sites ($n_s$, -) | Number of Routes ($n_R$, -) | Resource Budget ($b_r$, km) | LP Solution ($J_{LP}$, -) | Obtained Solution ($J_A$, -) | Optimality Gap (%) | Computing Time (s) | Generated Columns (-) |
|---|---|---|---|---|---|---|---|---|
| p1.4.e | 30 | 4 | 6.2 | 10.0 | 10 | 0.00 | 0.01 | 3 |
| p1.4.g | 30 | 4 | 8.8 | 35.0 | 35 | 0.00 | 0.01 | 8 |
| p1.4.i | 30 | 4 | 11.5 | 60.0 | 60 | 0.00 | 0.01 | 27 |
| p1.4.k | 30 | 4 | 13.8 | 95.0 | 95 | 0.00 | 0.02 | 91 |
| p1.4.m | 30 | 4 | 16.2 | 130.0 | 130 | 0.00 | 0.04 | 190 |
| p1.4.o | 30 | 4 | 18.2 | 155.0 | 155 | 0.00 | 0.06 | 336 |
| p1.4.q | 30 | 4 | 20.0 | 182.5 | 180 | 1.37 | 0.21 | 346 |
| p2.4.c | 19 | 4 | 5.8 | 70.0 | 70 | 0.00 | 0.01 | 18 |
| p2.4.e | 19 | 4 | 6.8 | 105.0 | 105 | 0.00 | 0.01 | 27 |
| p2.4.g | 19 | 4 | 8.0 | 120.0 | 120 | 0.00 | 0.15 | 58 |
| p2.4.i | 19 | 4 | 9.5 | 140.0 | 140 | 0.00 | 0.08 | 100 |
| p2.4.k | 19 | 4 | 11.2 | 180.0 | 180 | 0.00 | 0.01 | 113 |
| p3.4.e | 31 | 4 | 8.8 | 170.0 | 170 | 0.00 | 0.01 | 38 |
| p3.4.g | 31 | 4 | 11.2 | 210.0 | 210 | 0.00 | 0.02 | 157 |
| p3.4.i | 31 | 4 | 13.8 | 260.0 | 260 | 0.00 | 0.04 | 340 |
| p3.4.k | 31 | 4 | 16.2 | 320.0 | 320 | 0.00 | 0.08 | 343 |
| p3.4.m | 31 | 4 | 18.8 | 380.0 | 380 | 0.00 | 0.11 | 506 |
| p3.4.o | 31 | 4 | 21.2 | 480.0 | 480 | 0.00 | 0.35 | 623 |
| p3.4.q | 31 | 4 | 23.8 | 573.3 | 560 | 2.33 | 0.98 | 574 |
| p3.4.s | 31 | 4 | 26.2 | 610.0 | 610 | 0.00 | 1.45 | 1601 |
| p4.4.c | 98 | 4 | 17.5 | 255.9 | 251 | 1.93 | 0.61 | 725 |
| p4.4.e | 98 | 4 | 22.5 | 318.8 | 318 | 0.26 | 3.4 | 2067 |
| p4.4.g | 98 | 4 | 27.5 | 428.1 | 415 | 3.05 | 57 | 3095 |
| p4.4.i | 98 | 4 | 32.5 | 503.8 | 496 | 1.56 | 568 | 5323 |
| p4.4.k | 98 | 4 | 37.5 | 651.8 | 646 | 0.89 | 83,669 | 10241 |
| p5.4.c | 64 | 4 | 3.8 | 20.0 | 20 | 0.00 | 0.01 | 5 |
| p5.4.e | 64 | 4 | 6.2 | 70.0 | 70 | 0.00 | 0.01 | 19 |
| p5.4.g | 64 | 4 | 8.8 | 130.0 | 130 | 0.00 | 0.01 | 70 |
| p5.4.i | 64 | 4 | 11.2 | 205.0 | 205 | 0.00 | 0.03 | 283 |
| p5.4.k | 64 | 4 | 13.8 | 355.0 | 355 | 0.00 | 0.09 | 385 |
| p5.4.m | 64 | 4 | 16.2 | 520.0 | 520 | 0.00 | 1.38 | 919 |
| p5.4.o | 64 | 4 | 18.8 | 670.0 | 655 | 2.24 | 4.67 | 1750 |
| p5.4.q | 64 | 4 | 21.2 | 860.0 | 860 | 0.00 | 37.60 | 3041 |
| p6.4.c | 62 | 4 | 6.2 | 42.0 | 42 | 0.00 | 0.01 | 12 |
| p6.4.e | 62 | 4 | 8.8 | 114.0 | 114 | 0.00 | 0.21 | 75 |
| p6.4.g | 62 | 4 | 11.2 | 186.0 | 186 | 0.00 | 0.06 | 497 |
| p6.4.i | 62 | 4 | 13.8 | 330.0 | 330 | 0.00 | 0.33 | 1148 |
| p6.4.k | 62 | 4 | 16.2 | 522.0 | 522 | 0.00 | 0.96 | 2170 |
| p6.4.m | 62 | 4 | 18.8 | 714.9 | 708 | 0.96 | 9.70 | 5000 |
| p7.4.c | 100 | 4 | 15.0 | 32.0 | 32 | 0.00 | 0.01 | 4 |
| p7.4.e | 100 | 4 | 25.0 | 93.0 | 93 | 0.00 | 0.02 | 26 |
| p7.4.g | 100 | 4 | 35.0 | 190.0 | 190 | 0.00 | 0.10 | 181 |
| p7.4.i | 100 | 4 | 45.0 | 330.7 | 324 | 2.03 | 1.96 | 2068 |
| p7.4.k | 100 | 4 | 55.0 | 455.2 | 449 | 1.35 | 26 | 4148 |
| p7.4.m | 100 | 4 | 65.0 | 619.5 | 602 | 2.82 | 1,090 | 9215 |

# VI. Case Studies

This section introduces two realistic case studies involving real-world applications of the VRPVP; 1) routing for planetary surface exploration and 2) city tour design for a tourist group. Note that both applications consider the involvement of multiple stakeholders and the profit obtained by visiting a site is expressed as a vector, representing the different perspectives of the stakeholders.

## A. Case 1: Routing for Planetary Surface Exploration

The design of the routes for a rover exploring the surface of a planet (e.g. Mars) is selected as the first case study subject (Ahn et al. 2008; Lee and Ahn 2017). We can imagine various different stakeholder groups interested in the planetary surface exploration, and each of the groups may have its own objective associated with the surface exploration mission. For example, geologists would be interested in gathering soil samples, biologists would be searching for evidence of life, and certain investors might be looking for natural resources from the exploration.

**Table 2: Benefits of space programs, as perceived by different organizations (Bainbridge 2015)**

| Benefit category | Percent of responses in category | | | |
| --- | --- | --- | --- | --- |
|  | NESFA, % | CFF, % | AIAA, % | Combined, % |
| Technological | 15.6 | 19.9 | 40.5 | 28.9 |
| Scientific | 25.3 | 13.9 | 20.2 | 20.1 |
| Political | 12.4 | 15.5 | 16.5 | 14.9 |
| Economic | 10.0 | 13.5 | 14.4 | 13.0 |
| Psychological | 22.1 | 11.8 | 5.6 | 11.5 |
| Religious | 7.1 | 19.3 | 2.1 | 7.5 |
| Social | 7.6 | 6.1 | 1.3 | 4.1 |
| (Number of responses) | (340) | (296) | (620) | (1,256) |

Table 2 summarizes the results of a study on the benefit of spaceflight as perceived by members of different organizations; the New England Science Fiction Association (NESFA), the Committee for the Future (CFF) and the American Institute of Aeronautics and Astronautics (AIAA). It should be noted that there are many different categories of benefit from the planetary exploration, and the mixture of categories is diverse, depending on the organization.

For this planetary surface exploration case, four stakeholder groups are set up based on the aforementioned study results: 1) technological, 2) scientific, 3) political, and 4) economic groups. The profit values for the four profit categories are randomly chosen while the sums of profits for different categories are the same.

The locations of the exploration sites were determined based on a test problem presented by Chao et al. (1996). The cost (length) to complete an arc is defined as the Euclidean distance between the two sites associated with the arc; the effect of terrain was not considered. The time to complete the exploration at each site ($t_i$) was randomly selected between 0.5 and 2.0 hours. Table 3 and Table 4 summarizes the parameters used in this case study. The site locations, profit values (for different stakeholders) and mission completion time are presented Table 5.

**Table 3: Problem Instance Parameters for Case 1**

| Parameter | Numerical value |
|---|---|
| Number of Stakeholders (-) | 4 |
| Time required to obtain profit: $t_i$ (hour) | 0.5 ~ 2.0 (uniformly distributed) |
| Profit assigned per site (for each stakeholder) | 1 ~ 30 (uniformly distributed) |
| Number of sites: $n_S$ (-) | 100 |
| Maximum number of routes: $n_R$ (-) | 5 |

**Table 4: Resource Constraint Parameters for Case 1 (Lee and Ahn 2017)**

| Constraint Type | Resource Type | Resource Budget | Consumption Coefficient |
|---|---|---|---|
| On-Route | Time | $b_r$ = 10 hours | $c_d$ = 0.125 hour/km[7] |
| | | | $c_r$ = 1 hour/hour |
| In-Mission | Time | $b_m$ = 45 hours | $d_d$ = 0.125 hour/km |
| | | | $d_r$ = 1 hour/hour |

---

[7] This value is calculated as 1/(vehicle speed). The vehicle speed was determined by averaging the speeds of rovers used in Apollo missions (8 km/hour).

**Table 5: Site information for Case 1 – locations, profit values, and stay time**

| Site No. | Location | | Profit for each stakeholder ($p_i^k$) | | | | $t_i$, min. | Site No. | Location | | Profit for each stakeholder ($p_i^k$) | | | | $t_i$, min. |
|---|---|---|---|---|---|---|---|---|---|---|---|---|---|---|---|
| | x, km | y, km | k = 1 | k = 2 | k = 3 | k = 4 | | | x, km | y, km | k = 1 | k = 2 | k = 3 | k = 4 | |
| 0 | 35 | 35 | 0 | 0 | 0 | 0 | 0 | | | | | | | | |
| 1 | 15 | 30 | 26 | 9 | 2 | 12 | 14 | 51 | 53 | 43 | 14 | 9 | 9 | 9 | 12 |
| 2 | 55 | 5 | 29 | 16 | 23 | 16 | 20 | 52 | 57 | 48 | 23 | 3 | 20 | 9 | 7 |
| 3 | 31 | 52 | 27 | 15 | 10 | 6 | 16 | 53 | 15 | 47 | 16 | 36 | 10 | 13 | 18 |
| 4 | 60 | 12 | 31 | 14 | 9 | 13 | 7 | 54 | 14 | 37 | 11 | 16 | 12 | 5 | 12 |
| 5 | 8 | 56 | 27 | 8 | 19 | 5 | 19 | 55 | 26 | 35 | 15 | 27 | 19 | 13 | 14 |
| 6 | 13 | 52 | 36 | 2 | 11 | 23 | 13 | 56 | 18 | 24 | 22 | 8 | 13 | 18 | 20 |
| 7 | 6 | 68 | 30 | 1 | 7 | 9 | 17 | 57 | 25 | 24 | 20 | 10 | 20 | 23 | 6 |
| 8 | 21 | 24 | 28 | 11 | 23 | 8 | 9 | 58 | 22 | 27 | 11 | 23 | 41 | 16 | 11 |
| 9 | 56 | 39 | 36 | 9 | 8 | 6 | 19 | 59 | 25 | 21 | 12 | 9 | 19 | 29 | 19 |
| 10 | 55 | 54 | 26 | 5 | 1 | 2 | 7 | 60 | 18 | 18 | 17 | 13 | 23 | 23 | 6 |
| 11 | 16 | 22 | 41 | 22 | 27 | 3 | 11 | 61 | 41 | 49 | 10 | 41 | 7 | 16 | 12 |
| 12 | 4 | 18 | 35 | 12 | 18 | 21 | 13 | 62 | 35 | 17 | 7 | 12 | 18 | 9 | 18 |
| 13 | 28 | 18 | 26 | 31 | 25 | 17 | 18 | 63 | 25 | 30 | 3 | 18 | 8 | 3 | 6 |
| 14 | 26 | 27 | 27 | 36 | 13 | 20 | 17 | 64 | 20 | 50 | 5 | 1 | 27 | 19 | 15 |
| 15 | 55 | 45 | 13 | 5 | 3 | 5 | 19 | 65 | 10 | 43 | 9 | 5 | 16 | 10 | 19 |
| 16 | 55 | 20 | 19 | 26 | 27 | 18 | 12 | 66 | 30 | 5 | 8 | 14 | 17 | 7 | 20 |
| 17 | 55 | 60 | 16 | 3 | 5 | 8 | 9 | 67 | 5 | 30 | 2 | 21 | 19 | 26 | 11 |
| 18 | 30 | 60 | 16 | 11 | 36 | 27 | 18 | 68 | 45 | 65 | 9 | 6 | 15 | 10 | 9 |
| 19 | 20 | 65 | 12 | 9 | 6 | 12 | 9 | 69 | 65 | 35 | 3 | 5 | 11 | 35 | 20 |
| 20 | 50 | 35 | 19 | 25 | 7 | 16 | 18 | 70 | 65 | 20 | 6 | 10 | 9 | 28 | 19 |
| 21 | 30 | 25 | 23 | 12 | 2 | 30 | 9 | 71 | 64 | 42 | 9 | 28 | 16 | 16 | 20 |
| 22 | 15 | 10 | 20 | 23 | 5 | 3 | 9 | 72 | 63 | 65 | 8 | 11 | 15 | 5 | 12 |
| 23 | 10 | 20 | 19 | 26 | 28 | 15 | 8 | 73 | 2 | 60 | 5 | 13 | 17 | 25 | 17 |
| 24 | 20 | 40 | 12 | 3 | 9 | 7 | 16 | 74 | 20 | 20 | 8 | 30 | 22 | 9 | 15 |
| 25 | 15 | 60 | 17 | 19 | 16 | 11 | 20 | 75 | 40 | 25 | 9 | 17 | 5 | 26 | 13 |
| 26 | 45 | 20 | 11 | 6 | 9 | 16 | 11 | 76 | 42 | 7 | 5 | 2 | 8 | 6 | 16 |
| 27 | 45 | 10 | 18 | 7 | 12 | 2 | 8 | 77 | 24 | 12 | 5 | 16 | 6 | 8 | 18 |
| 28 | 45 | 30 | 17 | 10 | 25 | 9 | 9 | 78 | 23 | 3 | 7 | 3 | 2 | 16 | 7 |
| 29 | 35 | 40 | 16 | 6 | 3 | 10 | 11 | 79 | 2 | 48 | 1 | 13 | 30 | 7 | 13 |
| 30 | 41 | 37 | 16 | 7 | 26 | 13 | 12 | 80 | 49 | 58 | 10 | 9 | 10 | 9 | 9 |
| 31 | 40 | 60 | 21 | 18 | 18 | 20 | 13 | 81 | 27 | 43 | 9 | 19 | 6 | 25 | 14 |
| 32 | 35 | 69 | 23 | 11 | 9 | 18 | 12 | 82 | 63 | 23 | 2 | 2 | 5 | 27 | 13 |
| 33 | 53 | 52 | 11 | 19 | 35 | 14 | 11 | 83 | 53 | 12 | 6 | 14 | 16 | 14 | 12 |
| 34 | 65 | 55 | 14 | 16 | 21 | 22 | 10 | 84 | 32 | 12 | 7 | 13 | 6 | 18 | 16 |
| 35 | 5 | 5 | 16 | 15 | 18 | 15 | 19 | 85 | 17 | 34 | 3 | 20 | 25 | 8 | 9 |
| 36 | 11 | 14 | 18 | 20 | 31 | 17 | 20 | 86 | 27 | 69 | 10 | 17 | 5 | 1 | 14 |
| 37 | 6 | 38 | 16 | 16 | 3 | 36 | 10 | 87 | 15 | 77 | 9 | 9 | 14 | 3 | 7 |
| 38 | 47 | 47 | 13 | 27 | 9 | 6 | 20 | 88 | 37 | 47 | 6 | 17 | 9 | 36 | 12 |
| 39 | 37 | 31 | 14 | 16 | 26 | 31 | 19 | 89 | 37 | 56 | 5 | 9 | 3 | 19 | 6 |
| 40 | 57 | 29 | 18 | 5 | 10 | 18 | 14 | 90 | 44 | 17 | 9 | 26 | 17 | 27 | 10 |
| 41 | 36 | 26 | 18 | 7 | 7 | 14 | 10 | 91 | 46 | 13 | 8 | 6 | 13 | 1 | 14 |
| 42 | 12 | 24 | 13 | 7 | 12 | 11 | 17 | 92 | 61 | 52 | 3 | 27 | 18 | 9 | 7 |
| 43 | 24 | 58 | 19 | 8 | 20 | 20 | 19 | 93 | 56 | 37 | 6 | 18 | 14 | 5 | 9 |
| 44 | 62 | 77 | 20 | 35 | 11 | 19 | 14 | 94 | 11 | 31 | 7 | 16 | 11 | 10 | 10 |
| 45 | 49 | 73 | 25 | 29 | 29 | 26 | 9 | 95 | 26 | 52 | 9 | 3 | 16 | 11 | 11 |
| 46 | 67 | 5 | 25 | 9 | 26 | 2 | 16 | 96 | 31 | 67 | 3 | 18 | 16 | 25 | 18 |
| 47 | 57 | 68 | 15 | 10 | 14 | 11 | 6 | 97 | 15 | 19 | 1 | 18 | 13 | 19 | 15 |
| 48 | 47 | 16 | 25 | 25 | 16 | 41 | 13 | 98 | 22 | 22 | 2 | 25 | 9 | 7 | 9 |
| 49 | 49 | 11 | 18 | 20 | 36 | 12 | 19 | 99 | 19 | 21 | 10 | 8 | 3 | 17 | 9 |
| 50 | 49 | 42 | 13 | 19 | 8 | 3 | 20 | 100 | 20 | 26 | 9 | 23 | 1 | 9 | 17 |

Fig. 3 and Table 6 show the solution of VRPVP for Case 1. The square at the location (35, 35) represents the depot where the rover starts and ends each route (for refueling and maintenance) and the dots are exploration sites. The objective function values for the LP relaxation ($J_L$) and the approximate MILP ($J_A$) are respectively 325.59 and 318. The optimality gap is computed as 2.33 %, which is a very low value indicating that the solution can be considered nearly optimal.

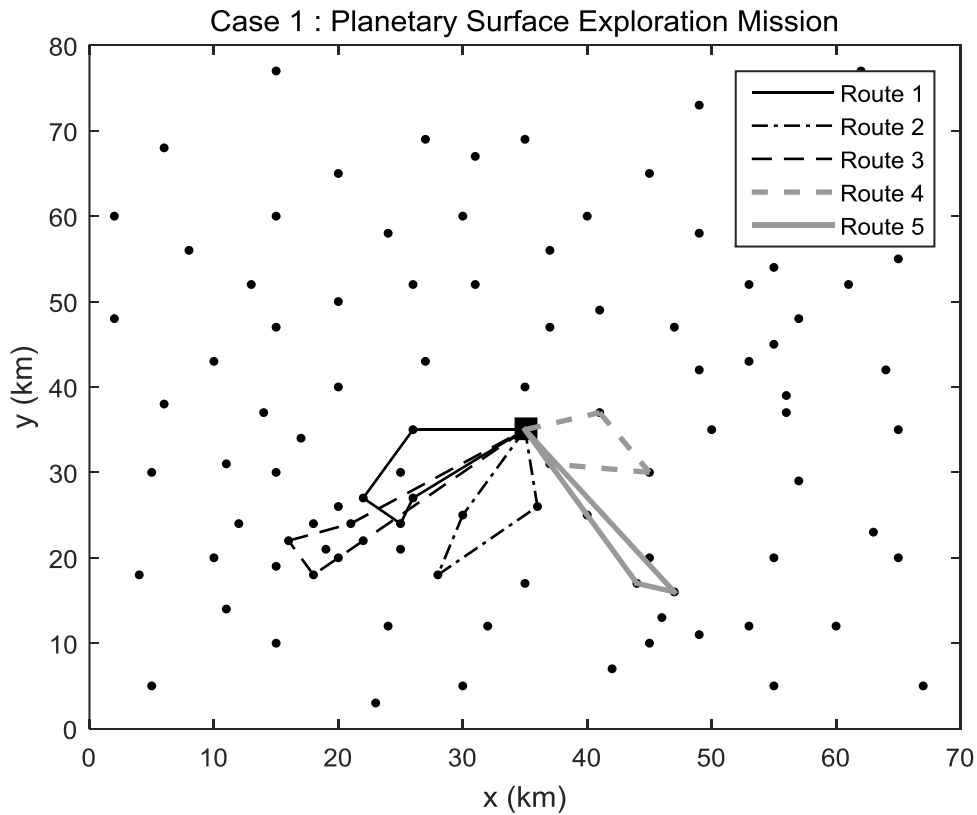

**Figure 3: VRPVP solution for Case 1 (planetary surface exploration case)**

**Table 6: VRPVP solution (routes) for Case 1**

| Route # | Visiting sites for routes | Travel time, hour | Collected profits (4 stakeholder groups) |
|---|---|---|---|
| Route 1 | Depot – 14 – 57 – 58 – 55 – Depot | 9.5 | 73 / 96 / 93 / 72 |
| Route 2 | Depot – 21 – 13 – 41 – Depot | 8.6 | 67 / 50 / 34 / 61 |
| Route 3 | Depot – 98 – 60 – 11 – 8 – Depot | 10.0 | 88 / 71 / 82 / 41 |
| Route 4 | Depot – 30 – 28 – 39 – Depot | 7.4 | 47 / 33 / 77 / 53 |
| Route 5 | Depot – 48 – 90 – 75 – Depot | 9.3 | 43 / 68 / 38 / 94 |
| | Total | 44.8 | **318** / **318** / 324 / 321 |

The profit sums for stakeholders obtained by solving the VRPVP were compared with the results of the VRPPs with various objective functions in Table 7. While each VRPP solution yields the maximum profits sum for the stakeholder associated with its objective function, the VRPVP provides the best minimum profit sum over all stakeholders.

**Table 7: Comparison of profits collected by VRPVP with VRPP results (Case 1)**

| Problem Type | Profit Sum for Stakeholder 1 | Profit Sum for Stakeholder 2 | Profit Sum for Stakeholder 3 | Profit Sum for Stakeholder 4 | Minimum Profit Sum | Sum of All Profit Sums |
|---|---|---|---|---|---|---|
| **VRPVP** | 318 | 318 | 324 | 321 | **318** | 1281 |
| **VRPP** | | | | | | |
| J = profit sum of stakeholder 1 | **349** | 182 | 234 | 207 | 182 | 972 |
| J = profit sum of stakeholder 2 | 231 | **379** | 225 | 285 | 225 | 1120 |
| J = profit sum of stakeholder 3 | 259 | 264 | **371** | 226 | 226 | 1120 |
| J = profit sum of stakeholder 4 | 254 | 264 | 242 | **401** | 242 | 1161 |
| J = sum of profit sums for all stakeholders | 315 | 334 | 313 | 340 | 313 | **1302** |

## B. Case 2: Design of Tourist Group Tour

The second case study solves a tour routing problem for a group of travelers with various interests using the VRPVP proposed in this paper. When people travel to a famous tourist spot as a group, they are faced with the problem of selecting the attractions to visit and determining their tour sequence. This routing problem influences the satisfaction level of the group significantly. A number of studies and their implementations that address this challenge can be found in the literature. For example, Vansteenwegen et al. (2007) proposed "the tourist trip design problem (TTDP)" and the mobile tourist guide, which can suggest holiday plans in real-time using reliable data on tourist attractions, was developed based on this problem.

One of the important advantages of their work is the capability to recommend a tour plan by considering the user's personal preferences on attractions. However, this approach is not appropriate for providing recommendations to a tourist group of multiple members. Since the preference structures of the group members could be highly diversified, the routes of the tour should be carefully determined so that the satisfaction levels of the members are harmonized. This tour routing problem, which considers the various preference structures of the tourist group members, was solved using the VRPVP framework introduced in this paper.

Tables 8 and Table 9 summarize the parameters used for the case study problem. Rome was selected as the location of the tour, which contains numerous attractions. In the case study, we selected 34 famous tourist attractions as candidate sites to visit. For each site/member, a profit value ($p_i^k$) is randomly assigned, while the sums of profit values for all members are identical. The time to be spent at each site ($t_i$) is determined based on the tourist guide website[8]. One of the five-star hotels located in the central area of the city was selected as the "depot." The location ($x_i$, $y_i$), stay time ($t_i$), and profit values for different members ($p_i^k$) associated with all sites are presented in Table 10. The interests of the travelers are reflected in assigning the profit values. We assumed that traveler 1 is interested in visiting religious places (e.g. churches), traveler 2 likes walkaway attractions (e.g. Palatine Hill), and traveler 3 is enthusiastic about ancient architectures (e.g. Foro Romano).

**Table 8: Problem Instance Parameters for Case 2**

| Parameter | Numerical value |
|---|---|
| Number of Stakeholders (-) | 3 |
| Time required to obtain profit: $t_i$ (hour) | 0.17 ~ 3.0 |
| Profit assigned per site (for each stakeholder) | 0 / 5 / 10 / 15 |
| Number of sites: $n_S$ (-) | 34 |
| Maximum number of routes: $n_R$ (-) | 3 |

**Table 9: Resource Constraint Parameters for Case 2**

| Constraint Type | Resource Type | Resource Budget | Consumption Coefficient |
|---|---|---|---|
| On-Route | Time | $\mathbf{b}_r$ = 5 hours | $\mathbf{c}_d$ = 1 hour/hour |
| | | | $\mathbf{c}_r$ = 1 hour/hour |

---

[8] The times required to obtain profit at each attraction are set based on the estimated time spent according to the Sygic Travel website (https://travel.sygic.com).

**Table 10: Site information for Case 2 –locations, profit values, and stay time**

| Site No. | Name of attraction | Location | | Profit for each stakeholder ($p_i^k$) | | | $t_i$, min. |
|---|---|---|---|---|---|---|---|
| | | latitude, deg. | longitude, deg. | Traveler 1 | Traveler 2 | Traveler 3 | |
| 0 | Hotel (Grand Hotel de la Minerve) | 41.8975 | 12.4777 | 0 | 0 | 0 | 0 |
| 1 | Colloseum | 41.8902 | 12.4922 | 1 | 1 | 2 | 180 |
| 2 | Arch of Constantine | 41.8898 | 12.4906 | 2 | 2 | 2 | 15 |
| 3 | Basilica of Saint Clement | 41.8893 | 12.4976 | 3 | 1 | 1 | 60 |
| 4 | Church of Four Crowned Martyrs | 41.8882 | 12.4990 | 3 | 1 | 2 | 15 |
| 5 | Tempio di Castore e Polluce | 41.8919 | 12.4857 | 1 | 1 | 3 | 60 |
| 6 | Foro Romano | 41.8925 | 12.4853 | 1 | 1 | 3 | 30 |
| 7 | Palatine Hill | 41.8895 | 12.4875 | 1 | 3 | 2 | 60 |
| 8 | House of Augustus | 41.8883 | 12.4867 | 1 | 1 | 1 | 60 |
| 9 | House of Livia | 41.8893 | 12.4856 | 1 | 1 | 1 | 60 |
| 10 | Forum of Caesar | 41.8939 | 12.4851 | 1 | 1 | 2 | 60 |
| 11 | Campidoglio Square | 41.8934 | 12.4828 | 2 | 2 | 2 | 15 |
| 12 | Santa Maria in Aracoeli | 41.8940 | 12.4832 | 1 | 1 | 1 | 30 |
| 13 | National Monument to Victor Emmanuel II | 41.8946 | 12.4831 | 1 | 1 | 1 | 30 |
| 14 | Trajan's Market | 41.8957 | 12.4862 | 1 | 3 | 3 | 150 |
| 15 | Trajan's Column | 41.8958 | 12.4843 | 1 | 2 | 2 | 15 |
| 16 | Venice Square | 41.8958 | 12.4826 | 1 | 2 | 1 | 30 |
| 17 | Church of Jesus | 41.8959 | 12.4799 | 3 | 1 | 1 | 15 |
| 18 | Theatre of Marcellus | 41.8919 | 12.4799 | 1 | 1 | 2 | 15 |
| 19 | Great Synagogue | 41.8920 | 12.4780 | 1 | 1 | 1 | 90 |
| 20 | Turtle Fountain | 41.8938 | 12.4776 | 1 | 2 | 2 | 15 |
| 21 | Mouth of Truth | 41.8881 | 12.4815 | 1 | 2 | 1 | 15 |
| 22 | Circus Maximus | 41.8861 | 12.4852 | 2 | 2 | 3 | 30 |
| 23 | Pantheon | 41.8986 | 12.4769 | 1 | 2 | 2 | 45 |
| 24 | Trevi Fountain | 41.9009 | 12.4833 | 1 | 3 | 1 | 45 |
| 25 | Capuchin Crypt | 41.9049 | 12.4884 | 3 | 1 | 1 | 60 |
| 26 | Pyramid of Cestius | 41.8765 | 12.4809 | 2 | 2 | 2 | 15 |
| 27 | Baths of Caracalla | 41.8790 | 12.4924 | 1 | 3 | 2 | 60 |
| 28 | Spanish Square & Spanish Steps | 41.9057 | 12.4823 | 2 | 3 | 1 | 30 |
| 29 | Basilica di San Pietro | 41.9022 | 12.4539 | 3 | 2 | 1 | 60 |
| 30 | Sistine Chapel | 41.9029 | 12.4544 | 3 | 2 | 1 | 60 |
| 31 | St Peters Square | 41.9022 | 12.4568 | 3 | 2 | 1 | 15 |
| 32 | Vatican Museums | 41.9066 | 12.4535 | 3 | 1 | 2 | 180 |
| 33 | Castel SantAngelo | 41.9031 | 12.4662 | 2 | 1 | 2 | 120 |
| 34 | Peoples Square | 41.9107 | 12.4764 | 1 | 1 | 1 | 15 |

The travel time between two sites, which is used as a cost associated with an arc connecting the two sites, was obtained using the *Google Maps Directions* API[9]. The TSP solution associated with site set $j$ (TSP$_j$) was expressed as the shortest time to complete the visits to all sites (in hours), not as the smallest path length.

Fig. 4 and Table 11 present the VRPVP solution for Case 2. The travel group visited 23 out of 34 tourist attractions using three routes that start/terminate at the depot (hotel). The obtained profit sum values for all tour group members were identical, with a value of 38. The numbers of attractions visited by the three routes are 6 (route 1), 7 (route 2), and 10 (route 3) while their travel times are all very close to the time budget for a route (5 hours). The values of the objective function for the LP relaxation ($J_L$) and the approximate MILP ($J_A$) are respectively 38.25 and 38. The optimality gap is computed as 0.66 %.

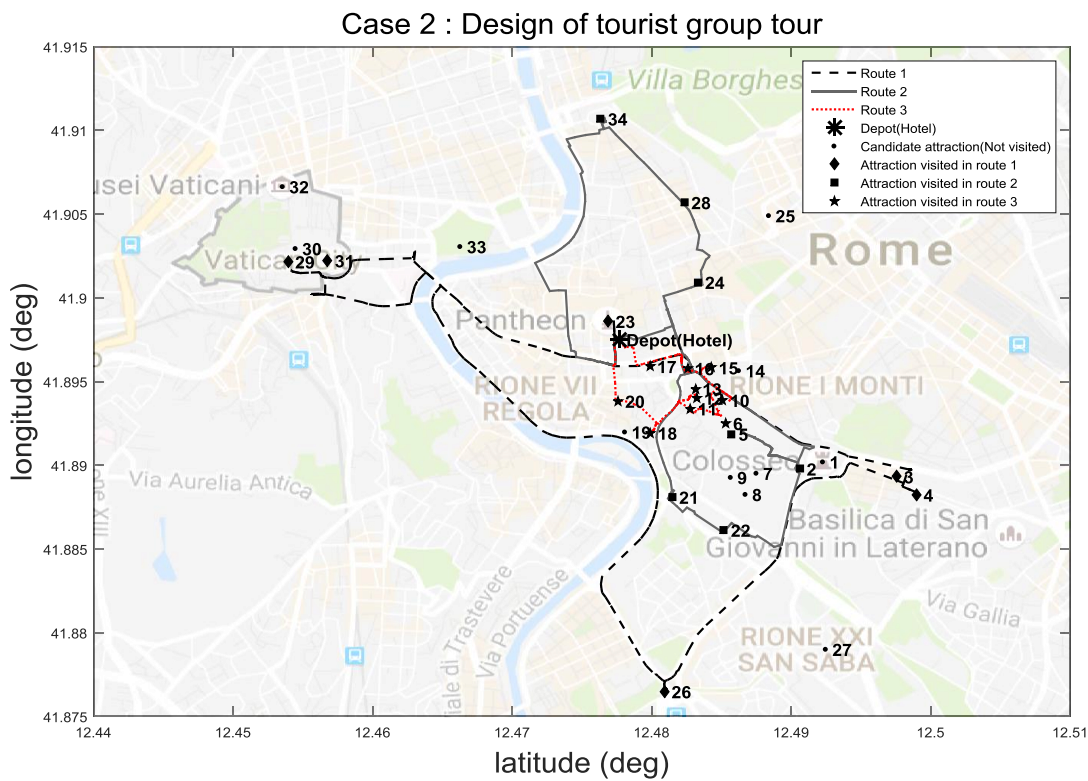

**Figure 4: VRPVP solution for Case 2 (Rome tour case)**

---

[9] API is an abbreviated form of "application programming interface."

Table 11: VRPVP solution (routes) for Case 2

| Route # | Visiting sites for routes | Travel time, hour | Collected profits (3 stakeholder groups) |
|---|---|---|---|
| Route 1 | Hotel – 31 – 29 – 26 – 3 – 4 – 23 – Hotel | 4.98 | 15 / 10 / 9 |
| Route 2 | Hotel – 5 – 2 – 22 – 21 – 24 – 28 – 34 – Hotel | 4.97 | 10 / 14 / 12 |
| Route 3 | Hotel – 17 – 16 – 15 – 10 – 6 – 11 – 13 – 12 – 18 – 20 - Hotel | 4.95 | 13 / 14 / 17 |
|  | Total | 14.9 | **38** / **38** / **38** |

The profit sums for stakeholders obtained by solving the VRPVP were compared with the results of the VRPPs with various objective functions in Table 12. These comparison results are similar to the comparison results for Case 1 presented in Table 7, which also demonstrates the effectiveness of the proposed framework in maximizing the profit sum of the least satisfied stakeholder.

Table 12: Comparison of profits collected by VRPVP with VRPP results (Case 2)

| Problem Type | Profit Sum for Member 1 | Profit Sum for Member 2 | Profit Sum for Member 3 | Minimum Profit Sum | Sum of All Profit Sums |
|---|---|---|---|---|---|
| **VRPVP** | 38 | 38 | 38 | **38** | 114 |
| **VRPP** | | | | | |
| J = profit sum for member 1 | **42** | 39 | 35 | 35 | 116 |
| J = profit sum for member 2 | 37 | **42** | 35 | 35 | 114 |
| J = profit sum for member 3 | 31 | 40 | **39** | 31 | 110 |
| J = sum of profit sums for all members | 40 | 41 | 36 | 36 | **117** |

# VII. Conclusions

The vehicle routing problem with vector profits (VRPVP), which is a routing problem that can handle the profit structures of multiple stakeholders associated with the mission, is proposed as an extension of the existing vehicle routing problem with profits (VRPP) framework. Maximizing the minimum profit sum (over all stakeholders) obtained through the routing was selected as the objective of the problem and the resource consumptions on individual routes and in the whole mission were considered as constraints. A solution method composed of the linear program (LP) relaxation and the column-generation technique that can generate the near-optimal solution of the VRPVP along with the optimality gap value was propose and validated through numerical experiments using test problem instances. Two case studies – the planetary surface exploration design and the routing for a group tour – demonstrated that the proposed VRPVP framework could be effectively used to solve routing problems that involve multiple stakeholders with satisfaction levels that should be reflected in their solutions in a balanced way.

The use of max-min criterion does not guarantee that the solution belongs to the Pareto front, which is one of issues that the following study should address. The authors expect that the issue can be resolved by introducing the lexicographic max-min criterion (Ogryczak 1997), which can be implemented by modifying the VRPVP framework proposed in this study. Another research direction is the development of a multi-objective optimization framework that can directly generate the family of efficient solutions (or, Pareto front) of the VRPVP. Extension of current framework to handle multiple depots and/or separate departure and arrival locations would be another area for potential future study.